\documentclass[10pt]{article}

\usepackage{amsmath,amstext,amsopn,amsfonts,eucal,amssymb,graphicx,multirow}

\newtheorem{theorem}{Theorem}[section]

\newtheorem{lemma}[theorem]{Lemma}
\newtheorem{proposition}[theorem]{Proposition}
\newtheorem{definition}[theorem]{Definition}

\begin{document}


\title{Counting Irreducible \\ Double Occurrence Words}

\author{Jonathan Burns\thanks{University of South Florida, \texttt{jtburns@mail.usf.edu}}
, Tilahun Muche\thanks{University of South Florida, \texttt{tmuche@mail.usf.edu}}}

\date{}

\maketitle


\begin{abstract}
A double occurrence word $w$ over a finite alphabet $\Sigma$ is a word  in which each 
alphabet letter appears exactly twice. Such words arise naturally in the study of topology, 
graph theory, and combinatorics. Recently, double occurrence words have been used for 
studying DNA recombination events. We develop formulas for counting and enumerating 
several elementary classes of double occurrence words such as palindromic, irreducible, 
and strongly-irreducible words.
\end{abstract}

\section{Introduction}

\begin{figure}[h]
\begin{center}$
\begin{array}{ccc}
\includegraphics[scale=.5]{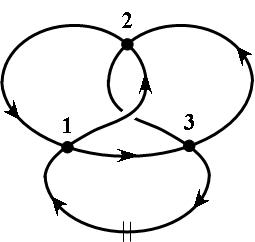} &
\includegraphics[scale=.5]{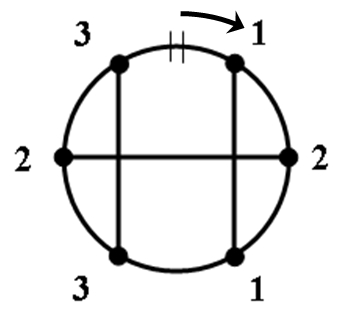} &
\includegraphics[scale=.5]{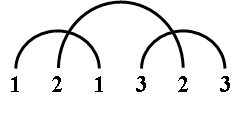} \\
\end{array}$
\end{center}
\label{diagrams}
\caption{Self-intersecting closed curve (left), chord diagram (center), and linked 
diagram (right) representations  of the double occurrence word $121323$. Base 
points, indicating the starting point for reading the word, are marked by $\parallel$.}
\end{figure}

A {\it double occurrence} word $w$ of size $n$ is a word containing $n$ distinct 
letters in any order which appear exactly twice, i.e., the {\it length} of $w$ is $2n$. 
There are three common pictorial representations of double occurrence words: 
self-intersecting closed curves in $\mathbb{R}^3$, chord diagrams, and linked 
diagrams as depicted in  Figure \ref{diagrams}.

Topologically, a double occurrence word with $n$ distinct letters can be interpreted as 
a closed curve traversing $n$ fixed points in $\mathbb{R}^3$ twice. Such a 
curve (also called an assembly graph \cite{Burns}) is self-intersecting and 
may contain over and under crossings when projected into the plane. Each curve 
of this type can be characterized through the double occurrence word corresponding 
to a  path following the direction of the curve in relation to a fixed {\it base point}. 
Self-intersecting closed curves are closely related to Gauss words, knot diagrams, and 
their shadows \cite{Cairns,Kauffman}.

Chord diagrams are defined in the following way. Start with a circle and place $n$ 
distinctly labeled chords with distinct endpoints in any arrangement (possibly crossing) 
around the circle. Label the endpoints of each chord with the chord  label. Fix a 
{\it base point} on the circle between any two chord endpoints on the circle. The 
resulting diagram is called a {\it chord diagram}. Each chord diagram has an associated 
double occurrence word formed by reading the labels of the endpoints, from the base 
point back to base point, clockwise around the circle. See \cite{Godsil,Klazar}  for more 
information on chord diagrams. 

A {\it linked} (or  {\it linearized chord} \cite{Stoimenow}) {\it diagram} is a pairing of $2n$ distinct 
ordered points. Graphically, the ordered points are positioned on a line and their
pairing is illustrated by an arc connecting them. Such a diagram can be specified 
by listing the pairs defined by the $n$ arcs. See \cite{Touchard1,Touchard2, Stein}. A linked diagram can be  obtained from a chord diagram by cutting the outer 
circle at the base point. Conversely, if we arrange the points of the  link diagram in
a circle and mark a base point between the first and last point, the corresponding 
representation is a chord diagram. 

Since double occurrence words naturally arise in a variety of contexts, insight into their 
combinatorial structure enriches several fields simultaneously. In this paper, we
explore several classifications of double occurrence words based on separating
larger double occurrence words into smaller double occurrence words. Further, we
count and enumerate members of these classes. 

Some of these formulas have been derived in completely different contexts using 
a variety of approaches. Moreover none of the papers we came across seemed 
to contain a compilation of the known formulas. In this paper we give a unified 
approach to deriving these formulas and provide a new formula, giving what 
appears to be an unobserved integer sequence.

We note that applications of double occurrence words extend to other 
disciplines. In \ref{separations}, we observe that certain double occurrence words
are related to particular Feynman diagrams in physics, and in Section \ref{biological} 
we establish a connection between double occurrence words and  DNA recombination 
events.

\section{Preliminaries}			

\subsection{Types of Equivalences}

For convenience, we let $\Sigma=\{1,2,\dots,n\}$ and relabel each double occurrence word such that when $i$ appears for the first time in the word, it is preceded by  $1,2,\dots,i-1$.
Double occurrence words labeled by this convention are said to be in {\it ascending order}. 
Two double occurrence words are said to be {\it equivalent} if they are equal after 
being relabeled in ascending order. If two double occurrence words are not equivalent,
they are said to be {\it distinct}. Throughout this paper, we shall assume 
that all double occurrence words are in ascending order unless stated otherwise.

For example, $122313$ is a double occurrence word in ascending order. Its {\it reverse} 
with the same letters is $313221$, which is not in ascending order. By relabeling  $313221$ in asscending order we obtain $121332$. In this example $122313$ is distinct from its reverse $121332$. However it is easily checked that $123312$ is equivalent to its reverse which motivates the following classification.

 \begin{definition} \label{palin_def}
{\rm
A double occurrence word is {\it palindromic} (or {\it symmetric}) if it is equivalent to its reverse.
A double occurrence word that is palindromic is called a {\it palindrome}. 
 }
 \end{definition}

In all three interpretations of double occurrence words (topological, graph theoretic, and linked diagrams), the reverse word induces a diagram, isomorphic to the original, with the orientation reversed. In the topological sense, the orientation refers to the orientation of the closed curve. While the reverse of a linked diagram may be interpreted as reading the diagram right-to-left rather than left-to-right. Finally, the reverse chord diagram may be achieved by reading the letters of the circle in a counter-clockwise fashion rather than clockwise.

If we wish to count the non-isomorphic diagrams generated from double occurrence words, we observe  that each diagram can have exactly two orientations. Thus, no more than two distinct double occurrence words can correspond to the same diagram with regard to a starting base point.

If a diagram corresponds to a palindrome, only one distinct double occurrence word is associated with the diagram. Therefore we may count the number of non-isomorphic diagrams with
regard to a base point as 
\begin{align*}
\mbox{Total Diagrams}\; & = (\mbox{\# of Palindromes})+\frac{1}{2}(\mbox{\#  of  Non-Palindromes}) \\[8pt]
 & = \frac{(\mbox{\# of D.O. Words}) + (\mbox{\# of Palindromes} )}{2}. \tag{$*$} \\
\end{align*}

\vspace{-.4cm} \noindent We will make use of this formula extensively throughout Section \ref{counting} to count the number of distinct diagrams corresponding to double occurrence words with each separation property.

It should be noted that omitting the base point in the closed curve or chord diagram makes it 
possible for more than two double occurrence words to be associated with the same diagram. For
instance, rotating the base point around the circle in Figure \ref{diagrams} would lead to 
121323, 213231, and 132312 which is 121323, 123132, and 123213 in ascending order,
respectively. We do not consider isomorphisms of this type in this paper.

\subsection{Types of Separations}
\label{separations}

As mentioned in the introduction, double occurrence words regularly appear in
various fields of mathematics. Unfortunately as a result, there are several different, 
and sometimes conflicting, definitions used to express identical properties. We shall make
note of these discrepancies in notation as they come up.

Jacques Touchard was one of the first researchers to comprehensively consider 
the counting of double occurrence words. In his paper \cite{Touchard1}, he classified 
several types of linked diagrams and enumerated the number of diagrams containing a 
fixed number of crossings. He introduced the classification of ``unique systems" and 
``proper unique systems" which coincide with the following two definitions for irreducible 
and strongly-irreducible words.

\begin{definition} \label{irr_def}
{\rm
If a double occurrence word $w$ can be written as a product $w=uv$ of two
non-empty double occurrence words $u, v$, then $w$ is called {\it 
reducible}; otherwise, it is called {\it irreducible}.
}
\end{definition}

The number of irreducible double occurrence words has a close connection with the number
of non-isomorphic unlabeled connected Feynman diagrams (also called {\it irreducible Feynman diagrams} \cite{Robinson}) arising in a simplified model of quantum electrodynamics \cite{Cyitanovic,Martin}.

This definition for irreducibility agrees with \cite{Angeleska} and \cite{Burns} yet conflicts with \cite{Stein} where ``irreducible"  is used for our notion of strongly-irreducible as defined below.

\begin{definition} \label{strong_irr_def}
{\rm A non-empty double occurrence word is {\it strongly-irreducible} if it does not contain a proper sub-word that is also a double occurrence word.}
 \end{definition}

The double occurrence word $12213434$ is reducible because it can be written as the 
product of the two double occurrence words $1221$ and $3434$, but $12344123$ is 
irreducible. However, since $44$ is a proper sub-word of $12344123$ it is not 
strongly-irreducible. The word $12132434$ is strongly-irreducible. By definition, 
strongly-irreducible words are also irreducible, so $12132434$ is irreducible as well. 
In particular 11 is strongly-irreducible.

Strongly-irreducible double occurrence words are also called {\it connected} words 
\cite{Klazar}. This terminology is motivated by the {\it circle graph} associated with a 
chord diagram. The circle graph is formed by representing the chords as vertices and 
the intersection of those chords as edges in the graph. In the topological convention, 
a circle graph is also called an {\it interlinking graph} \cite{Cairns}. Without too many 
difficulties it can be proven that a double occurrence word is strongly-irreducible if and 
only if the circle graph of the corresponding chord diagram, or interlinking graph of the
corresponding closed curve, is connected.

\begin{lemma}
\label{sub-word}
{\rm Every double occurrence word contains a strongly-irreducible sub-word.}
\end{lemma}
{\it Proof.} If a double occurrence word $w$ is strongly-irreducible, then $w$ itself 
is a strongly-irreducible sub-word of $w$. Double occurrence words which are not 
strongly-irreducible, by definition, contain a proper sub-word $w_1$ which is a double 
occurrence word and is either strongly-irreducible or not. If the sub-word is not 
strongly-irreducible we check the reducibility of its proper sub-word $w_2$. Since $w$ 
has finite length, we must reach a double occurrence word $w_i$, which is a 
strongly-irreducible proper sub-word of $w_{i-1}$, through finite recursion. Since 
$w_i$ must be a proper sub-word of $w$, this completes the proof. \hfill $\Box$

\section{Counting}
\label{counting}

It is well known \cite{Burns,Klazar,Stein,Touchard2} and straightforward 
to show that the total number of double occurrence words is $(2n-1)!!$. Formula 
$(*)$ motivates us to enumerate the number of double occurrence words which 
correspond to palindromes.

\subsection{Palindromes}

\begin{theorem}
The number $L_n$ of palindromic double occurrence words of \\ length $2n$, is given by
$$L_n = \sum\limits_{k=0}^{\lfloor n/2 \rfloor} \frac{n!}{(n-2k)! \, k!} \quad \mbox{for }n\ge1. $$
\end{theorem}
{\it Proof.\/} 
Observe that $L_1 = 1$ since there is a unique one letter palindrome, and $L_2 = 3$ 
because 1122, 1212, and 1221 are all the two letter palindromes.

If a double occurrence word $w$ of size $n \ge 2$ is a palindrome beginning 
and ending with 1, then the word formed by removing both $1$s is also a 
palindrome. Hence there are $L_{n-1}$ palindromes with $n$ letters that 
start and end with 1.

Now consider a word $w$ of size $n \ge 3$ where the second symbol $1$ is 
at the position $j \ne 2n$. Note that there are $2n-2$ possible positions for $j$.
Then the word $w$ is a palindrome if and only if $w$ contains the same symbol 
$s$ at the positions $2n$ and $n - j + 1$. Removing symbols $1$ and $s$ from 
$w$, and relabeling the resulting word accordingly, produces a palindrome of 
length $n - 2$. Hence there are $L_{n-2}$ palindromes that have a 
symbol $1$ at the $j$th position for $2 \le j \le 2n-1$.

According to the above argument, $$L_n = L_{n-1} + (2n - 2) L_{n-2} 
\mbox{ for } n \ge 3, \quad L_1 = 1 \mbox{ and } L_2 = 3$$ is a recurrence relation 
for $L_n$. It is known~\cite{OEIS} that the closed formula for this recursive relations 
is as stated. \hfill $\Box$

This formula is expressed without proof in a comment by Ross Drewe in $A047974$ of the 
OEIS \cite{OEIS} in 2008, but this may not be the original source. Similar results, such as 
the number of palindromic chord diagrams without a base point, were known in 2000 \cite{Stoimenow}. The above proof reprinted here is found in \cite{Burns}.

\subsection{Irreducibles}
Though Touchard introduced the classification of irreducible words in 1952, there 
seems to be little continuation of his efforts. In 2000, Martin and Kearney \cite{Martin} 
expressed the number of irreducible words in the broader context of solutions to 
generating functions. Here, we address the count and construction of both the 
irreducible double occurrence words and irreducible palindromes directly.

\begin{lemma}
The number of irreducible double occurrence words $I_n$ with length $2n$ satisfies 
the recurrence formula  $I_1 = 1$ and $$ I_n = (2n-1)!! - \sum\limits_{k=1}^{n-1} 
I_{n-k} \;(2k-1)!! \quad \mbox{for } n \ge 2.$$
\end{lemma}
{\it Proof.\/} We shall count the number of irreducible double occurrence words by 
subtracting the number of reducible double occurrence words from the total 
number of double occurrence words of length $2n$ and show that each reducible 
word may be written as the product of an irreducible word and a non-empty 
double occurrence word. 

Without loss of generality, let $w=uv$ be a reducible double occurrence word
of length $2n$ such that $u$ is also an irreducible double occurrence word. 
Note that every proper prefix of an irreducible word is not necessarily a double 
occurrence word. If the length of $v$ is $2k$, for some $1 \le k \le n-1$, then 
the length of $u$ is $2(n-k)$. By construction, $u$ is irreducible and is counted 
among $I_{n-k}$ and $v$ is counted among the $(2k-1)!!$ possible double 
occurrence words of length $2k$.

Summing over the possible symbols in $v$ yields the desired count. Since $u$ 
is irreducible and $v$ is non-empty, this ensures that each reducible double 
occurrence word $w$ is counted exactly once. \hfill $\Box$

\begin{theorem}
The number of irreducible palindromes $J_n$ with length $2n$ satisfies the recurrence formula  $J_1 = 1$ and $$ J_n = L_n - \sum\limits_{k=1}^{\lfloor n/2 \rfloor} (2k-1)!! 
\; J_{n-2k} \quad \mbox{for } n \ge 2.$$ where $L_n$ is the total number of 
palindromes with length $2n$.
\end{theorem}
{\it Proof.\/} Similar to the above argument, we first count the reducible 
palindromes and subtract them from the total number of palindromic words.

Suppose $w$ is a reducible double occurrence word with length $2n$. Then 
$w$ can be written as $w=uvu'$ where $u$ is an arbitrary double occurrence 
word with length $2k$ ($1 \le k \le \lfloor n/2 \rfloor$), $u'$ is the double 
occurrence word corresponding to $u$ by reversing the orientation, and $v$ 
is an irreducible palindrome  with length $2(n-2k)$. \hfill $\Box$

Though the number of irreducible double occurrence words appears in the 
OEIS (A000698), we note that the number of irreducible palindromes is the 
only sequence discussed in this paper which is not currently listed in the OEIS  
\cite{OEIS}. See Table \ref{do-words}, Table \ref{pal-words}, and Table 
\ref{iso-classes} for the number of irreducibles, strong-irreducibles, and the number of non-isomorphic diagrams as defined according to ($*$), respectively.

\subsection{Strong-Irreducibles}
The classification of strongly-irreducible double occurrence words was introduced in 
\cite{Touchard2} and the first counting of the strong-irreducibles was done by Stein 
in \cite{Stein}. Stein was the first to count both the strongly-irreducible double occurrence
words and the strongly-irreducible palindromes, but his counting methods and recursive
formulas were simplified in \cite{Nijenhuis} and later by Klazar in \cite{Klazar}. In Theorem
\ref{strongly_irr_theorem}, we present a proof  similar to \cite{Klazar} expressed 
in terms of language theory.

Using language theory to count double occurrence words led directly to a characterization
of the strongly-irreducible double occurrence words, which we express in Lemma
\ref{strong_irr_lemma}, and Theorem \ref{strongly_irr_theorem} follows as a natural 
consequence.

\begin{lemma}\label {strong_irr_lemma}
Every strongly-irreducible double occurrence word $w$ in ascending order may be written 
in a unique form as $w=1u_{1}v_{1}1v_{2}u_{2}$ where $1 u_1 1 u_2$ and $v_1 v_2$ 
are both strongly-irreducible.
\end{lemma}
{\it Proof}. Let $w$ be strongly-irreducible. Every double occurrence word $w$ in ascending order must be of the form $w=1p_{1}1p_{2}$. Delete both 1's. Then we have a double occurrence word  $p_{1}p_{2}=u_{1}xu_{2}$ where $x$ is the  first strongly-irreducible double occurrence word of smallest positive length. Thus $u_{1}$ and $u_{2}$ are uniquely defined. Note that $u_{1}$ and $u_{2}$  may be empty words.

Let $v_1$ be the prefix of $x$ which is a suffix of $p_1$ and let $v_2$ be the suffix of $x$ 
which is the prefix of $p_2$. This means that $x=v_1 v_2$. Neither $v_1$ nor $v_2$ is 
empty as it would imply that $x$ is a sub-word of either $p_1$ or $p_2$ which would constitute
a proper sub-word of $w$. Since $w$ is taken to be strongly-irreducible, this cannot be.

We show that $1u_{1}1u_{2}$ is strongly-irreducible. Suppose not. Then there exists a 
non-empty double occurrence sub-word ${\it z}$ in either $u_{1}$ or $u_{2}$ which 
implies that $w$ contains $z$ and is not strongly-irreducible. This is a contradiction. Hence 
$1u_{1}1u_{2}$ and $v_1 v_2$ are strongly-irreducible. \hfill $\Box$

\begin{theorem} 
\label{strongly_irr_theorem}
The number of strongly-irreducible double occurrence words $S_{n}$ with 
length $2n$ satisfies the recurrence formula  $$S_{n} =(n-1) \sum_{k=1}^{n-1}S_{k}S_{n-k},$$ where 
$S_{1}=1$ and $n\geq 2.$
\end{theorem}
{\it Proof.} 
Note that the only strongly-irreducible double occurrence word of length 2 is $11$, 
i.e., $S_1= 1$. 

Let $u$ and $v$ be strongly-irreducible double occurrence words such that the 
length of $v$ is $2k$, the length of $u$ is $2(n-k)$, $u =1 u_1 1 u_2$,  
and $v=v_1 v_2$. Since the length of $v$ is $2k$, there are $2k-1$ ways to write 
$v=v_{1}v_{2}$ with $v_{1},v_{2}$  not empty.  By Lemma \ref{strong_irr_lemma}, 
each strongly-irreducible double occurrence word $w$ of length $2n$ can be uniquely 
represented as $w= 1u_{1}v_{1}1v_{2}u_{2}$. Hence there are $2k-1$ possibilities 
for such $w$'s to be formed from each $u$ and $v$. 

Since there are $S_{n-k}$ choices for $u$ and $S_k$ choices for $v$ the total 
counting for $S_n$ when $n \geq 2$ is given by $$S_{n}=\sum_{k=1}^{n-1}(2k-1)S_{k}S_{n-k}
= (n-1) \sum\limits_{k=1}^{n-1} S_{k}S_{n-k}. \eqno{\Box}$$

\noindent For completeness, we state Klazar's counting formula of the strongly-irreducible palindromes. See \cite{Klazar} for the proof.

\begin{theorem}
\label{strongly_irr_pal_theorem}
Let $S_n$ and $T_n$ be the number of strongly-irreducible double occurrence words and strongly-irreducible palindromes of length $2n$, respectively. Then $$T_n = \sum\limits_{i=1}^{n-2} T_i T_{n-i}+\sum\limits_{i=1}^{\lfloor n/2 \rfloor}(2n-4i-1) S_i T_{n-2i}$$
for $n \ge 2$ where $T_0 = -1$ and $T_1 = 1$.
\end{theorem}

Theorem \ref{strongly_irr_theorem} and Theorem \ref{strongly_irr_pal_theorem} correspond to
the sequences A000699 and A004300 listed in the OEIS. For the first few values of these sequences, see Table \ref{do-words} and Table \ref{pal-words}.

\section{Connection with DNA recombination}
\label{biological}

Several species of ciliates, such as {\it Oxytricha} and {\it Stylonychia}, undergo 
massive genome rearrangement during sexual reproduction. These massively occurring  recombination processes make them ideal model organisms to study gene rearrangements. 
See \cite{Ehrenfeucht} and references therein for details of the descriptions below.

There are two types of nuclei, a micronucleus and a macronucleus, in these species. 
Micronuclear genes contain both coding and non-coding segments which are reassembled to 
macronuclear genes during sexual reproduction. The coding segments, called
{\it macronuclear destined sequences} or {\it MDSs}, are part of  the final 
unscrambled gene. The individual MDSs within a micronuclear gene may be separated
by non-coding segments, called {\it internal eliminated sequences} or {\it IESs}, which 
are excised during the recombination process.

\begin{figure}
\centering
\includegraphics[scale=.4]{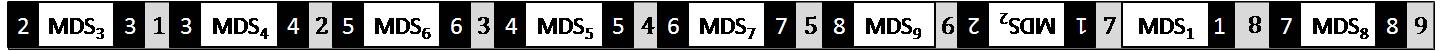}
\caption{Scrambled Actin I micronuclear gene in {\it Oxytricha nova} \cite{Prescott2}.}
\label{gene1}
\end{figure}

\begin{figure}
\centering
\includegraphics[scale=.4]{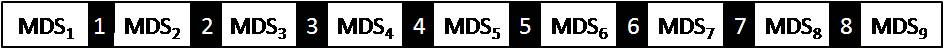}
\caption{Unscrambled Actin I macronuclear gene in {\it Oxytricha nova} \cite{Prescott2}.}
\label{gene2}
\end{figure}

In relation to an unscrambled macronuclear gene (Fig. \ref{gene2}), a scrambled 
micronuclear gene (Fig. \ref{gene1}) may have permuted or inverted MDS segments 
separated by IESs. Formation of the macronuclear genes in these ciliates thus requires 
any combination of the following three events: unscrambling of segment order, DNA 
inversion, and IES removal.

Since the IESs are removed in the unscrambled gene, it is only necessary to record the
order and direction of the MDSs in the scrambled gene. A {\it micronuclear arrangement} (cf. \cite{Ehrenfeucht}) is a sequence of permuted and inverted MDSs. In particular, each micronuclear arrangement $\alpha$ with $k$ MDSs has a corresponding permutation $\sigma_{\alpha}:[k] \rightarrow [k]$ and a signing function $\epsilon_\alpha:[k] \rightarrow \{-1,+1\}$ which uniquely defines the arrangement. A sign of $-1$ indicates that an MDS is inverted with respect to the gene sequence in the macronuclear gene while a sign of $+1$ indicates a regular orientation.

For example, the micronuclear arrangement of the Actin I gene in Figure \ref{gene1} is
$$M_3^{+1} M_4^{+1} M_6^{+1} M_5^{+1} M_7^{+1} M_9^{+1} M_2^{-1} M_1^{+1} M_8^{+1}$$ or more commonly denoted $$M_3 M_4 M_6 M_5 M_7 M_9 \overline{M}_2 M_1 M_8$$ where $\overline{M}_2$ indicates that $\mbox{MDS}_2$ is inverted in the scrambled micronuclear gene.

\begin{proposition}
Let $A_n$ be the number of micronuclear arrangements of $n$ MDSs. Then $$A_n = 2^{n}n! = (2n)!!.$$
\end{proposition}
{\it Proof.} Each micronuclear arrangement $\alpha$ with $n$ MDSs is uniquely 
defined by its corresponding permutation $\sigma_\alpha$ and signing function 
$\epsilon_\alpha$. Since each MDS may be signed in one of two ways, there are 
$2^{n}$ ways to sign the $n!$ permutations of all arrangements of $\alpha$ with 
$n$ MDSs. \hfill $\Box$

\bigskip

The exact process by which the scrambled micronuclear gene recombines into an unscrambled macronuclear gene is unknown. However it is theorized \cite{Prescott} that short sequences 
of nucleotides, called {\it pointers}, found at the beginning and end of each MDS, guide the recombination process. In fact, each MDS is characterized by its pointers in the following sense.

Each MDS is labeled according to its order in the unscrambled macronuclear gene. The 
pointers flanking the MDSs correspond to the order of the MDSs such that the pointer 
sequence at the end of the  $i$th MDS coincides with the pointer sequence at the 
beginning of the $(i+1)$th MDS. IESs are excised and their coding is not necessary. 
Since the pointers at the beginning and end of the whole gene do not align with any 
other pointers, we omit them. Mathematically, this translates to the following.

Let $\mathcal{A}_n$ be the set of all micronuclear arrangements with $n$ MDSs and $\mathcal{K}_n$ be the set of all double occurrence words with length $2n$. Then $\varrho:\mathcal{A}_n \rightarrow \mathcal{K}_{n-1}$ is a homomorphism which translates a micronuclear arrangement to the ordered sequence of pointers which describes it, i.e.,

\begin{enumerate}
\item $\varrho(M_1^{-1}) \mapsto (1)$ and $\varrho(M_1^{+1}) \mapsto (1)$
\item $\varrho(M_i^{-1}) \mapsto (i)(i-1)$
\item $\varrho(M_i^{+1}) \mapsto (i-1)(i)$
\item  $\varrho(M_n^{-1}) \mapsto (n-1)$ and $\varrho(M_n^{+1}) \mapsto (n-1).$
\end{enumerate}
For the micronuclear arrangement $\alpha = M_2^{-1}M_4^{+1}M_1^{+1}M_5^{-1}M_3^{+1}$, $$\varrho(\alpha) = (2)(1)(3)(4)(1)(4)(2)(3)$$ which corresponds to the double occurrence 
word $12342413$ in ascending order. Therefore each scrambled micronuclear gene 
corresponds to a  micronuclear arrangement which, in turn, has an associated double
occurrence word. 

A double occurrence word is called {\it realizable} if it has a corresponding micronuclear
arrangement. The shortest double occurrence word which is not realizable is $11233244$.
For further information on realizable double occurrence words see \cite{Burns}.

\section{Conclusions}

Double occurrence words are studied in topology, graph theory, and combinatorics 
by way of self-intersecting closed curves in $\mathbb{R}^3$,  chord graphs and linked 
diagrams, respectively. Their applications extend beyond abstraction to other 
disciplines such as physics and genetics. We considered the counting and enumeration of 
several reducibility classes of double occurrence words which directly led to a new 
characterization of strongly-irreducible double occurrence words. Further, all but one of the enumerated sequences are listed in the OEIS \cite{OEIS}, which suggests both the 
relevance of the previously listed enumerations and the novelty of the unlisted irreducible palindrome count. It should be noted that all the counting arguments present in this 
paper followed a similar theme: separate the classes of double occurrence words into 
palindromes and non-palindromes and describe the construction of large double 
occurrence words from smaller double occurrence words. We believe that the counting
techniques presented here could be used to enumerate new classes of double occurrence 
words as they arise in future research.

\section{Acknowledgments}
This research was supported through the efforts of the NSF Grant DMS  \#0900671.

We would like to thank our research associates Egor Dolzhenko, Nata\v{s}a Jonoska, 
and Masahico Saito for their assistance in improving both the mathematical and 
written content of this paper.

\begin{table}[p]
\centering
\begin{tabular}{|c r r r |}
\hline
Symbols & All	& Irreducible	& Strongly Irreducible	\\
\hline
\hline
1	&	1		&		1		&	1		\\
2	&	3		&		2		&	1		\\
3	&	15		&		10		&	4		\\
4	&	105		&		74		&	27		\\
5	&	945		&		706		&	248		\\
6	&	10395		&		8162		&	2830		\\
7	&	135135	&		110410	&	38232		\\
8	&	2027025	&		1708394	&	593859	\\	
9	&	34459425	&		29752066	&	10401712	\\	
10	&	654729075	&		576037442	&	202601898	\\	
11	&	13749310575	&		12277827850	&	4342263000	\\
12	&	316234143225	&	285764591114	&	101551822350	\\
\hline \hline
OEIS	&	A001147 $(K_n)$	&	A000698 $(I_n)$	&	A000699 $(S_n)$	\\
\hline
\end{tabular}
	\caption{All Double Occurrence Words.}\label{do-words}
\end{table}

\begin{table}[p]
\centering
\begin{tabular}{|c r r r |}
\hline
Symbols & All	& Irreducible	& Strongly Irreducible	\\
\hline
\hline
1	& 1				&	1				&	1 \\
2	& 3				&	2				&	1 \\
3	& 7				&	6				&	2 \\
4	& 25			&	20				&	7 \\
5	& 81			&	72				&	22 \\
6	& 331			&	290				&	96 \\
7	& 1303			&	1198			&	380 \\
8	& 5937			&	5452			&	1853 \\
9	& 26785			&	25176			&	8510 \\
10	& 133651		&	125874			&	44940 \\
11	& 669351		&	637926			&	229836 \\
12	& 3609673		&	3448708			&	1296410 \\
\hline \hline
OEIS	&	A047974 $(L_n)$	&	------ $(J_n)$	&	A004300 $(T_n)$	\\
\hline
\end{tabular}
	\caption{Palindromic Double Occurrence Words.}\label{pal-words}
\end{table}

\begin{table}[h]
\centering
\begin{tabular}{|c r r r |}
\hline
Symbols 	&	 All	&	Irreducible	&	Strongly Irreducible	\\
\hline
\hline
1	&	1	&	1	&	1       	\\
2	&	3	&	2 	&	1	\\
3	&	11	&	8	&	3       	\\
4	&	65	&	47	&	17	\\
5	&	513	&	389	&	135	\\
6	&	5363	&	4226	&	1463	\\
\hline
\hline
\multirow{2}{*}{OEIS}	&	A001147 $(K_n)$ 	&	A000698 $(I_n)$	&	A000699 $(S_n)$	\\
	&	A047974 $(L_n)$	&	------ $(J_n)$	&	A004300 $(T_n)$	\\
\hline
\end{tabular}
	\caption{Non-isomorphic diagrams in ($*$) are obtained by summing all words 
	with the palindromes of each class and halving the total. These sequences
	do not appear in the OEIS \cite{OEIS}, but can be built from listed 
	sequences.
} 
\label{iso-classes}
\end{table}

\pagebreak

\pagebreak

\end{document}